%% Plain TeX
%%

\newcount\secno
\newcount\prmno
\newif\ifnotfound
\newif\iffound

% Les commandes
\def\section#1{\vskip0,8truecm
               \global\def\currenvir{section}
               \global\advance\secno by1\global\prmno=0
               {\bf \number\secno. {#1}}
               \smallskip}

\def\subsection{\global\def\currenvir{subsection}
                \global\advance\prmno by1
       \smallskip           \ind{ (\number\secno.\number\prmno) }}
\def\subsec{\global\def\currenvir{subsection}
                \global\advance\prmno by1
                { (\number\secno.\number\prmno)\ }}

\def\proclaim#1{\global\advance\prmno by 1
                {\bf #1 \the\secno.\the\prmno$.-$ }}

\long\def\th#1 \enonce#2\endth{%
   \medbreak\proclaim{#1}{\it #2}\global\def\currenvir{th}\smallskip}

\def\bib#1{\rm #1}
\long\def\thr#1\bib#2\enonce#3\endth{%
\medbreak{\global\advance\prmno by 1\bf#1\the\secno.\the\prmno\ 
\bib{#2}$\!\!.-$ } {\it
#3}\global\def\currenvir{th}\smallskip}
%usage:\thr Proposition \bib Dupont \enonce
\def\rem#1{\global\advance\prmno by 1
{\it #1} \the\secno.\the\prmno$.-$ }

\magnification 1250
\pretolerance=500 \tolerance=1000  \brokenpenalty=5000
\mathcode`A="7041 \mathcode`B="7042 \mathcode`C="7043
\mathcode`D="7044 \mathcode`E="7045 \mathcode`F="7046
\mathcode`G="7047 \mathcode`H="7048 \mathcode`I="7049
\mathcode`J="704A \mathcode`K="704B \mathcode`L="704C
\mathcode`M="704D \mathcode`N="704E \mathcode`O="704F
\mathcode`P="7050 \mathcode`Q="7051 \mathcode`R="7052
\mathcode`S="7053 \mathcode`T="7054 \mathcode`U="7055
\mathcode`V="7056 \mathcode`W="7057 \mathcode`X="7058
\mathcode`Y="7059 \mathcode`Z="705A
\def\spacedmath#1{\def\packedmath##1${\bgroup\mathsurround =0pt##1
\egroup$}
\mathsurround#1
\everymath={\packedmath}\everydisplay={\mathsurround=0pt}}
\def\nospacedmath{\mathsurround=0pt
\everymath={}\everydisplay={} } \spacedmath{2pt}
\def\qfl#1{\buildrel {#1}\over {\longrightarrow}}
\def\phfl#1#2{\normalbaselines{\baselineskip=0pt
\lineskip=10truept\lineskiplimit=1truept}\nospacedmath\smash 
{\mathop{\hbox to 8truemm{\rightarrowfill}}
\limits^{\scriptstyle#1}_{\scriptstyle#2}}}
\def\hfl#1#2{\normalbaselines{\baselineskip=0truept
\lineskip=10truept\lineskiplimit=1truept}\nospacedmath\smash
{\mathop{\hbox to
12truemm{\rightarrowfill}}\limits^{\scriptstyle#1}_{\scriptstyle#2}}}
\def\diagram#1{\def\normalbaselines{\baselineskip=0truept
\lineskip=10truept\lineskiplimit=1truept}   \matrix{#1}}
\def\vfl#1#2{\llap{$\scriptstyle#1$}\left\downarrow\vbox to
6truemm{}\right.\rlap{$\scriptstyle#2$}}

\def\iso{\vbox{\hbox to .8cm{\hfill{$\scriptstyle\sim$}\hfill}
\nointerlineskip\hbox to .8cm{{\hfill$\longrightarrow $\hfill}} }}
\def\pprod_#1^#2{\raise 2pt
\hbox{$\mathrel{\scriptstyle\mathop{\kern0pt\prod}\limits_{#1}^{#2}}$}}

\newfam\bboardfam

\font\eightrm=cmr8         \font\eighti=cmmi8
\font\eightsy=cmsy8        \font\eightbf=cmbx8 
\font\eighttt=cmtt8        \font\eightit=cmti8
\font\eightsl=cmsl8        \font\sixrm=cmr6   
\font\sixi=cmmi6           \font\sixsy=cmsy6
\font\sixbf=cmbx6\catcode`\@=11
\def\eightpoint{%
  \textfont0=\eightrm \scriptfont0=\sixrm \scriptscriptfont0=\fiverm
  \def\rm{\fam\z@\eightrm}%
  \textfont1=\eighti  \scriptfont1=\sixi  \scriptscriptfont1=\fivei
  \def\oldstyle{\fam\@ne\eighti}\let\old=\oldstyle
  \textfont2=\eightsy \scriptfont2=\sixsy \scriptscriptfont2=\fivesy
  \textfont\itfam=\eightit
  \def\it{\fam\itfam\eightit}%
  \textfont\slfam=\eightsl
  \def\sl{\fam\slfam\eightsl}%
  \textfont\bffam=\eightbf \scriptfont\bffam=\sixbf 
  \scriptscriptfont\bffam=\fivebf
  \def\bf{\fam\bffam\eightbf}%
  \textfont\ttfam=\eighttt
  \def\tt{\fam\ttfam\eighttt}%
  \abovedisplayskip=9pt plus 3pt minus 9pt
  \belowdisplayskip=\abovedisplayskip
  \abovedisplayshortskip=0pt plus 3pt
  \belowdisplayshortskip=3pt plus 3pt 
  \smallskipamount=2pt plus 1pt minus 1pt
  \medskipamount=4pt plus 2pt minus 1pt
  \bigskipamount=9pt plus 3pt minus 3pt
  \normalbaselineskip=9pt
  \setbox\strutbox=\hbox{\vrule height7pt depth2pt width0pt}%
  \normalbaselines\rm}\catcode`\@=12

\newcount\noteno
\noteno=0
\def\up#1{\raise 1ex\hbox{\sevenrm#1}}
\def\note#1{\global\advance\noteno by1
\footnote{\parindent0.4cm\up{\number\noteno}\
}{\vtop{\eightpoint\baselineskip12pt\hsize15.5truecm\noindent
#1}}\parindent 0cm}

\def\pc#1{\tenrm#1\sevenrm}
\def\tx{\kern-1.5pt -}
\def\cqfd{\kern 2truemm\unskip\penalty 500\vrule height 4pt depth 0pt 
width 4pt\medbreak} 
\def\virg{\raise
.4ex\hbox{,}}

\def\ind{\par\hskip 0.8truecm\relax}
\def\indp{\par\hskip 0.5truecm\relax}
\def\smallsetminus{\mathrel{\hbox{\vrule height 3pt depth -2pt width 6pt}}}

\def\rond{\kern 1pt{\scriptstyle\circ}\kern 1pt}
\def\det{\mathop{\rm det}\nolimits}
\def\aut{\mathop{\rm Aut}\nolimits}
\def\Pic{\mathop{\rm Pic}\nolimits}
\def\rk{\mathop{\rm rk\,}\nolimits}
\def\Sy{\mathop{\rm Sym}\nolimits}
 
\frenchspacing
\input amssym.def
\input amssym
\vsize = 25truecm
\hsize = 16truecm
%\hoffset = -.15truecm
\voffset = -.5truecm
\parindent=0cm
\baselineskip15pt
\overfullrule=0pt
\def\P{\Bbb P}

\def\C{\Bbb C}
\def\Z{\Bbb Z}
\def\A{\Bbb A}

\def\pg{PGL_2(K)}

\def\au{\aut(\P^1_K)}
\def\auc{\aut(\P^1_{k(t)})}
\font\gragrec=cmmib10 
\def\mug{\hbox{\gragrec \char22}}
\def\ca{\mathop{\rm Cr}_k}
\font\grit=cmmib10
\def\pgrit{\grit\char112}

\centerline{\bf {\pgrit}\hskip1pt-elementary subgroups of the
Cremona group}

\smallskip
\smallskip \centerline{Arnaud {\pc BEAUVILLE}} 
\vskip1.2cm

{\bf Introduction}
\smallskip
\ind Let $k$ be an algebraically closed field. The {\it Cremona group}
$\ca$ is the group of birational transformations of $\P^2_{k}$, or
equivalently the group of $k$\tx automorphisms of the field
$k(x,y)$. There is 
an extensive classical literature about this group, in particular
about its finite subgroups -- see the introduction
of [dF] for a list of references.
\ind The classification of conjugacy classes of elements
of prime order $p$ in $\ca$ has been given a modern treatment
in [B-B] for $p=2$ and in [dF] for $p\ge 3$ (see also [B-Bl]). In
this note we go one step further and classify $p$\tx elementary
subgroups -- that is, subgroups isomorphic to 
$(\Z/p)^r$ for $p$ prime. We will mostly describe such a
subgroup as a group $G$ of automorphisms of a rational surface $S$:
we identify $G$ to a subgroup of
$\ca$ by choosing a birational map
$\varphi:S\dasharrow \P^2$. Then the conjugacy class of $G$ in $\ca$
depends only on the data $(G,S)$.

\smallskip  {\bf Theorem}$.-$ {\it
Let
$G$ be a subgroup of
$\ca$ of the form $(\Z/p)^r$ with $p$ prime $\not= {\rm char}(k)$.
Then:
\indp {\rm a)} Assume $p\ge 5$. Then $r\le 2$, and if $r=2$ $G$ is
conjugate to the $p$\tx torsion subgroup of the diagonal torus\note{By the {\it
diagonal torus} of $PGL_r(k)$ we mean the subgroup of
projective transformations $(X_0,\ldots ,X_r)\mapsto (t_0X_0,\ldots
,t _rX_r)$ for $t _0,\ldots ,t _r$ in $k^*$.}
 of $PGL_3(k)=\aut(\P^2)$.

\indp {\rm b)} Assume $p=3$. Then $r\le 3$, and if $r=3$ $G$ is
conjugate to the $3$\tx torsion subgroup of the diagonal torus
 of $PGL_4(k)$, acting on the Fermat cubic surface $X_0^3+\ldots
+X_3^3=0$.
\indp {\rm c)} Assume $p=2$. Then $r\le 4$, and if equality holds $G$
is conjugate to one of the following subgroups:
\ind ${\rm c}_1)$ the subgroup of $\ca$ spanned by the involutions:
$$(x,y)\mapsto (-x,y)\quad,\quad (x,y)\mapsto ({1\over
x},y)\quad,\quad(x,y)\mapsto (x,-y)\quad,\quad(x,y)\mapsto (x,{P(x)\over
y})$$
with $P(x)=\pprod_{a\in I}^{}(x^2+x^{-2}-a)$ for some finite subset $I$ of
$k\smallsetminus\{2,-2\}$.
\ind ${\rm c}_2)$
the $2$\tx torsion subgroup of the diagonal torus
 of  $PGL_5(k)$, acting on the quartic del Pezzo surface in
$\P^4$ with equations
$\sum_{i=0}^4X_i^2=\sum_{i=0}^4\lambda_iX_i^2=0$ for some distinct
elements $\lambda_0,\ldots \lambda_4$ of} $k$.
\smallskip\ind Note that the question does not make sense when
$p={\rm char}(k)$,  already in dimension 1: the group $PGL_2(k)$
contains a subgroup isomorphic to $k$, hence  infinite-dimensional over
$\Z/p$.
\ind 
In the next section we  discuss some
motivation for this question; then we reduce the problem  through
standard techniques to the study of $p$\tx elementary subgroups
$G\i\aut(S)$, where $S$ is either a del Pezzo surface or carries a
$\P^1$\tx fibration preserved by $G$. We will study the latter case in
\S 2 and the former in \S 3. Finally in \S 4 we discuss the
classification of the conjugacy classes of  subgroups isomorphic to
$(\Z/2)^4$   (the only case where the conjugacy class is not unique).
\ind As I. Dolgachev pointed out to us, the result (over $\C$) could be
deduced from the list of the finite subgroups of the Cremona group 
established by Kantor [K], and completed by Wiman [W]. However, the
results of \S 2 and \S 4   would still be needed to decide whether certain
subgroups are conjugate or not. Most of the results of \S 3 are contained
in those of Kantor and Wiman, but they are so much simpler in our 
specific situation that we have preferred to give  an independent proof.

 {\bigskip {\eightrm\baselineskip=12pt
\leftskip1cm\rightskip1cm\hskip0.8truecm  I am  indebted to
J.-P. Serre for asking me the question considered here, and for
useful discussions, in particular about (1.1) and (1.2)
below. I thank the referee for pointing out a gap in the
first version of this note.\par}}\smallskip  

\section{Comments, and beginning of the proof}
\subsection Though
very large,   the Cremona group  behaves in some respect like a
semi-simple group of rank 2:  every maximal torus  has
dimension 2, and is conjugate to the diagonal torus $T$ of
$PGL_3(k)$  [D1].  In this 
set-up the analogue of the Weyl group   is the
whole automorphisms group
$GL_2(\Z)$ of $T$ ([D1], cor. 5 p. 522). 
\ind  Now let
$H$ be a semi-simple  group over
$k$, and  $p$ a prime number
$\ge 7$ which does not divide the order of $\pi _1(H)$. Then every 
maximal $p$\tx elementary subgroup  $G\i H$ is the $p$\tx torsion
subgroup of a maximal torus  of $H$; moreover this torus  is
unique, in fact it is the centralizer of  $G$ in $H$ [Bo].
\subsection Our theorem shows that the first part of this statement also 
holds for the Cremona group (with $p\ge 5$). However, {\it the
maximal torus containing $G$ is  not unique}. Indeed, let $T$ be the
 diagonal torus of
$PGL_3(k)$, and $G$ its $p$\tx torsion subgroup.
The centralizer of $G$ in $\ca$ contains the transformations
$\sigma_f :(x,y)\mapsto (x,yf(x^p))$ for  $f\in k(t)^*$. 
If $f$ is not a monomial
 $\sigma _f$ does not normalize $T$, and  $G$ is
also contained in $\sigma_f T\sigma_f ^{-1} $. 
\smallskip 
\subsection
Now let us begin the proof of the Theorem. Let $G$ be a finite subgroup
of $\ca$. Then $G$ can be realized as a group of automorphisms of a
rational surface $S$ (see for instance [dF-E], Thm. 1.4). Moreover we can
assume that
$(G,S)$ is minimal, that is,  every birational $G$\tx equivariant
morphism of $S$ onto a smooth surface with a $G$\tx  action  is
an isomorphism.  Then one of the following holds:
\indp $\bullet$ $G$ preserves a fibration $f:S\rightarrow \P^1$ with
rational fibers;
\indp $\bullet$ $\rk\Pic(S)^G=1$.
\ind This result goes back to Manin [M], at least in the case (of interest
for us) when $G$ is abelian. It is by now a direct consequence
of Mori theory, see for instance [Z], lemma 4.3.
\ind In the former case $G$ embeds in the group of automorphisms of
the generic fibre $\P^1_{k(t)}$ of $f$; in the next section we are going
to classify the $p$\tx elementary subgroups of $\aut(\P^1_{k(t)})$. In
the latter case $S$ is a del Pezzo surface, and the group $\aut(S)$ is
well-known; we will use this information in \S 3 to classify the
corresponding $p$\tx elementary subgroups. Putting these
results together gives the Theorem.

\section{Subgroups of ${\bf Aut(P^1_K)}$}
\ind  Let $K$ be an extension of $k$, and $p$ a prime number
$\not= {\rm char}(k)$. Let us recall the classification of $p$\tx
elementary subgroups of $\pg$. Let  $C_p\i
PGL_2(k)$ be the cyclic subgroup of homographies
$z\mapsto \zeta z$, with $\zeta^p =1$; for $\delta \in K^*$, let
$V_\delta\i \pg $ be the subgroup of order 4 spanned by the
homographies $z\mapsto -z$ and $z\mapsto {\delta \over z}$. 

\th Lemma
\enonce  Let $G$ be a subgroup of
$\pg$ of the form $(\Z/p)^r$, with $p$ prime $\not= {\rm char}(k)$.
Then  $G$ is conjugate to a subgroup of $C_p$ if $p$ is odd,
and to a subgroup of $V_\delta$ for some $\delta \in K^*$ if $p=2$.
\endth
{\it Proof} : Let $\sigma$ be an element of $G $, represented by a matrix $A\in
GL_2(K)$ which satisfies $A^p=\delta I$ for some scalar $\delta \in K^*$. 
  Taking determinants give $(\det
A)^p=\delta ^2$, so that $\delta ^2\equiv 1$ (mod. $K^{*p}$).  If $p$ is
odd, this implies $\delta \in K^{*p}$; thus   $A$ is diagonalizable, and
$\sigma
$ is conjugate to an element of $C_p$. The centralizer of $C_p$ in $\pg$ is
the group of homotheties $z\mapsto \lambda z$, hence $G$ is contained
in its $p$\tx torsion subgroup  $C_p$.

\ind We now consider the case $p=2$. In a basis  of
$K^2$ of the form $(v,Av)$, we have $A=\pmatrix{0 & \delta \cr 1&0}$,
so that $\sigma $ becomes the homography $z\mapsto {\delta \over z}$. 
The centralizer of $\sigma $ in  $PGL_2(\overline{K})$ is
isomorphic to $\Z/2\times \overline{K}^*$, so its 2-torsion subgroup
 is isomorphic to $(\Z/2)^2$; since this subgroup contains $V_\delta $,
it is equal to $V_\delta$, and therefore 
 $G\i V_\delta $.\cqfd
\ind The conjugacy class of the subgroup $V_\delta $ depends only on the
class of $\delta $ in $K^*/K^{*2}$. In particular, when $K=k$, all these
subgroups are conjugate to $V_1$.
 
\subsection We are interested in the
automorphism group
$\au$ of the
$k$\tx scheme $\P^1_K$. Let $\Gamma $ be the automorphism group 
of $K$ over $k$. The
 action of $\au $ on the global functions of
$\P^1_K$ gives rise to an exact sequence
$$1\rightarrow \pg\longrightarrow \au \qfl{\pi} \Gamma 
\rightarrow 1\ ;$$ the surjection $\pi $ has a canonical section
$s:\Gamma
\rightarrow
\au$ which maps $\sigma \in \Gamma $ to the automorphism
$z\mapsto \sigma (z)$. 
We will use this section to identify $\Gamma $ to
a subgroup of $\au$, and view $\au$ as the semi-direct
product $\pg\rtimes \Gamma $. In other words, an element
of $\au$ is a transformation $\displaystyle z\mapsto {a\sigma (z)+b\over
c\sigma (z)+d}$ for some $\pmatrix{a&b\cr c&d}\in GL_2(K)$ and $\sigma
\in \Gamma $.
\smallskip 
\ind We are interested in the case $K=k(t)$, so that $\Gamma =PGL_2(k)$. 

\th Proposition
\enonce
Let $G$ be a subgroup of
$\auc$ of the form $(\Z/p)^r$, with $p$ prime $\not= {\rm char}(k)$.
\ind {\rm a)} When $p$ is odd, we have
$r\le 2$; if $r=2$, $G$ is conjugate to the
subgroup $C_p\times C_p$ of $\pg\rtimes PGL_2(k)$.
\ind {\rm b)} When $p=2$ we have $r\le 4$; if equality holds, $G$ is conjugate
to the subgroup $V_\delta \times V_1$ of $\pg\rtimes PGL_2(k)$ for some
element $\delta $ of $K^*$ invariant under $V_1$.
\endth

{\it Proof} : 
 Let $G'=G\cap \pg$. The bound on $r$ follows 
from the exact sequence
$$1\rightarrow G'\longrightarrow G\qfl{\pi }\pi (G)\rightarrow 1$$
and Lemma 2.1. Assume that  $p$ is odd and $r=2$. Then $G'$ and
$\pi (G)$ are cyclic, and after conjugation 
we may assume $G'=C_p$ and $\pi (G)=C_p$. 
 Let $\sigma $ be a generator of $C_p$, and $g$ an element of $G$ with 
$\pi (g)=\sigma $. 
We have $g=h\sigma $ for some element $h$ of $\pg$. Since the elements of
$G'$ commute with $\sigma $ and $h\sigma $, they must commute
with $h$; thus $h$ belongs to the centralizer of $G'=C_p$ in
$\pg$, hence is a homothety
$z\mapsto\lambda z$ for some $\lambda\in K^*$. 
 The relation $(h\sigma
)^p=1$ gives $\lambda \cdot
\sigma (\lambda )\cdot \ldots \cdot \sigma ^{p-1}(\lambda )=1$; by
Hilbert's theorem 90,
$\lambda $ is then equal to $\mu^{-1}  \sigma (\mu ) $ for some  $\mu
\in K^*$. The homothety $z\mapsto \mu z$ conjugates $g $ to
$\sigma $, and commutes with $G'$, hence conjugates $G$ to
$C_p\times C_p$.\smallskip \ind Now we assume $p=2$ and $r=4$. After 
conjugation  we may assume  $\pi (G)=V_1$, and $G'=V_\zeta $ for some
$\zeta \in K^*$ (Lemma 2.1). Let $\sigma \in V_1$. The elements of $\pi
^{-1} (\sigma )$ are of the form $h\sigma $ for some $h$ in $\pg$. Since $h$
must commute with the homography $z\mapsto -z$, it is of the form
$z\mapsto \lambda z$ or $z\mapsto {\mu\over z}$. Imposing that $h\sigma $
commutes with the homography $z\mapsto {\zeta  \over z}$ gives $\pi
^{-1}(\sigma )=G'\cdot g_\sigma $, where $g_\sigma (z)=\lambda_\sigma
\sigma (z)$, and $\lambda_\sigma \in K^*$ is a square root of
$\zeta \,\sigma (\zeta)^{-1} $. Note that $\lambda_\sigma $ is well-defined
up to sign. 
\ind We can choose 
each $\lambda_\sigma $ so that $\sigma \mapsto g_\sigma $ is a section of
$\pi $ (fix $\lambda_\sigma $ arbitrarily on a basis of $V_1$ over $\Z/2$, and
extend  linearly). Then $\sigma \mapsto\lambda_\sigma $ is a 1-cocycle of
$V_1$ with values in $K^*$. As above this cocycle is a coboundary, hence there
exists
$\mu\in K^*$ such that $\lambda_\sigma=\mu^{-1}\sigma (\mu)  $ for
each $\sigma \in V_1$. The homothety $z\mapsto \mu z$ conjugates
$g_\sigma $ to $\sigma $ and $V_\zeta $ to $V_{\delta}$, with $\delta
=\zeta\mu^2$. The commutativity of $G$ imposes that $\delta $ is invariant
under
$V_1$; this implies in particular that  $V_1 $ acts trivially on
$V_\delta $, hence
$G=V_\delta\times V_1\i\pg\rtimes PGL_2(k)$.\cqfd
\subsection Thus for $p$ odd and $r=2$, $G$ is conjugate to the  group of
diagonal transformations
$(z,t)\mapsto (\alpha z,\beta  t)$, with $\alpha,\beta \in\mug_p$. For $p=2$
and $r=4$, $G$ is conjugate to the group $V_\delta \times V_1$ generated by
$$(z,t)\mapsto (z,-t)\quad,\quad (z,t)\mapsto (z,\,{1\over
t})\quad,\quad(z,t)\mapsto (-z,t)\quad,\quad(z,t)\mapsto ({\delta\over
z}\,,t)\ .$$
 Since the map $t\mapsto t^2+t^{-2}$ identifies $\P^1/V_1$ with
$\P^1$, the invariance of $\delta $ under $V_1$ means that $\delta $ can be
written $Q(t^2+t^{-2})$ for some $Q\in k(s)$. Replacing $Q$ by $QF^2$ does
not change the conjugacy class of $V_\delta \times V_1$, so we can assume
that $Q$ is a polynomial with simple roots. Moreover conjugation by
 the birational transformation  $(z,t)\mapsto (z(t\pm  t^{-1}),t )$ 
amounts to multiply $Q$ by $(s\pm 2)$, so we can
assume that $Q$ does not vanish at $\pm 2$.
This gives case ${\rm
c}_1)$ of the  Theorem.

\section{Automorphisms of del Pezzo surfaces}
\ind We now consider the case where $S$ is a del Pezzo surface and
$G\cong (\Z/p)^r$ a subgroup of $\aut(S)$ such that $\rk\Pic(S)^G=1$. 
 We first recall  the following well-known fact, which is a particular case
of the results mentioned in (1.1):
\th Lemma
\enonce Let $G$
be a subgroup of $PGL_n(k)$ of the form $(\Z/p)^r$, with $p$ prime
$\not= {\rm char}(k)$. Assume that $p$ does not divide $n$. Then $r\le
n-1$, and if equality holds $G$ is conjugate to the $p$\tx torsion
subgroup of the diagonal torus.
\endth 
{\it Proof} : Pulling back $G$ to
$SL_n(k)$ gives a central extension of $G$ by the group $\mug_n$ of
$n$\tx th roots of unity in $k$.
Such extensions are parametrized by the group
$H^2(G,\mug_n)$ which is annihilated both by $n$ and $p$. Thus our
extension is trivial, and $G$ lifts to a subgroup  of
$SL_n(k)$, isomorphic to $(\Z/p)^r$. Such a subgroup is contained in a
maximal torus of $SL_n(k)$, hence our assertions.\cqfd
\subsection Let us start with the case $S=\P^2$. By the above lemma,
if $p\not= 3$, we have $r\le 2$, and any subgroup of $PGL_3(k)$
isomorphic to $(\Z/p)^2$ is conjugate to the $p$\tx torsion
subgroup of the diagonal torus.
\ind The case $p=3$ is classical (see e.g. [Bo], 6.4 for a more
general statement): we have again $r\le 2$, and a subgroup 
isomorphic to $(\Z/3)^2$ is conjugate either to the diagonal
subgroup, or to the subgroup spanned by the automorphisms
$(X_0,X_1,X_2)\mapsto (X_1,X_2,X_0)$ and $(X_0,X_1,X_2)\mapsto
(X_0,\alpha X_1,\alpha^2 X_2)$ for
$\alpha\in\mug_3$.\smallskip 
\subsection If $S$ is obtained from $\P^2$ by blowing up one or two points,
the group $\aut(S)$ is a subgroup of $PGL_3(k)$, so (3.2) applies.
Suppose that $S$ is the blow up of $\P^2$ at 
three  non-collinear points.
 Let $N \subset PGL_3(k)$ be the subgroup of  
automorphisms preserving this three-points subset.  The group
$\aut(S)$ is the semi-direct product of 
$N$ and a subgroup of order 2. Let $G$  be a subgroup of
$\aut(S)$ isomorphic to $(\Z/p)^r$. If $p\not= 2$, $G$ is contained in
$N$, so that (3.2) applies. If $p=2$, we have $r\le 3$ again by
(3.2).

\smallskip 
\ind We now consider the case $S=\P^1 \times \P^1$. 
\th Lemma
\enonce Let $G$ be a group of automorphisms of $\P^1\times \P^1$,
isomorphic to $(\Z/p)^r$, such that $\rk\Pic(\P^1\times \P^1)^G=1$.
Then $p=2$ and $r\le 3$.
\endth
{\it Proof} : The automorphism group of $\P^1\times \P^1$ is the
semi-direct product \break $\bigl(PGL_2(k)\times
PGL_2(k)\bigr)\rtimes
\Z/2$, where $\Z/2$ acts on $PGL_2(k)\times PGL_2(k)$ by
exchan\-ging the factors. If $p\not= 2$, our subgroup $G$ is contained
in 
$PGL_2(k)\times PGL_2(k)$, hence $\Pic(\P^1\times
\P^1)^G=\Pic(\P^1\times \P^1)$ has rank 2. 
\ind Thus we have $p=2$. The subgroup $G' $of  $G$ preserving the two
$\P^1$\tx fibrations is contained in $PGL_2(k)\times PGL_2(k)$, thus
in $C_2\times C_2$ up to conjugacy, and $G$ is conjugate to a
subgroup of the semi-direct product $(C_2\times C_2)\rtimes \Z/2$.
But the elements of order 2 in this group not contained in $C_2\times
C_2$ are contained in the 
(direct) product of $\Z/2$ by the diagonal subgroup $C_2\i C_2\times
C_2$. Therefore $G$ must be contained in $C_2\times \Z/2$, hence
$r\le 3$.\cqfd 

\subsection   It remains to consider the case when $S$ is 
obtained from $\P^2$ by blowing up $\ell $ points in general position,
with $4\le \ell \le 8$.
 We start by recalling some classical facts about such
surfaces, which can be found for instance in [D2]. The 
primitive cohomology $H^2(S,\Z)_{prim}$ (the orthogonal of the 
canonical class in $H^2(S,\Z)$) is the root lattice of a root system $R$.
The group $\aut(S)$ acts faithfully on
$H^2(S,\Z)_{prim}$, hence can be identified with a subgroup of the
automorphism group of the root system $R$; it is actually contained in
the Weyl group $W\i\aut(R)$ [Do]. 
The root systems which appear are the following (see [B]):
\indp $\bullet$ $\ell =4$: $R=A_4$, $W={\goth S}_5$
\indp $\bullet$ $\ell =5$: $R=D_5$, $W=(\Z/2)^4\rtimes {\goth S}_5$
\indp $\bullet$ $\ell =6$: $R=E_6$, $|W|=2^7.3^4.5$
\indp $\bullet$ $\ell =7$: $R=E_7$, $|W|=2^{10}.3^4.5.7$
\indp $\bullet$ $\ell =8$: $R=E_8$, $|W|=2^{14}.3^5.5^2.7$
\subsection For $\ell =7$ (resp. $\ell =8$), the linear system $|-K_S|$ 
(resp.
$|-2K_S|$) defines a degree 2 morphism onto $\P^2$ (resp. a quadric
cone in $\P^3$), branched along a canonically embedded smooth curve
$C$ of genus $3$ (resp. 4). This gives rise to an exact
sequence
$$1\rightarrow \Z/2\longrightarrow \aut(S)\longrightarrow
\aut(C)\eqno{(3.6)}$$ (the right hand side map is actually surjective, but
we will not need this). To control subgroups of $\aut(C)$ the following 
lemma will be useful:
\th Lemma
\enonce Let $p$ be a prime number, $r$ an integer, and $C$  a curve of
genus $g$ with a faithful action of the group $(\Z/p)^r$. Then
$p^{r-1}$ divides
$2g-2$. 
\endth
{\it Proof} : Put $G=(\Z/p)^r$, and consider the covering 
$\pi :C\rightarrow C/G$. Since the stabilizer of any point of $C$ is
cyclic, the fibre of $\pi $ above a branch point consists of $p^{r-1}$
points with ramification index $p$. Thus Hurwitz's
formula gives the result.\cqfd

\th Proposition
\enonce For $p$ prime $\ge 5$, the group $(\Z/p)^2$
does not act faithfully on a del Pezzo surface.
\endth
{\it Proof} :  A glance at the list 3.5 shows that the only case
where the order of the Weyl group is divisible by
$p^2$, with $p\ge 5$, is $\ell =8$. In this case the exact sequence
(3.6) shows that $(\Z/p)^2$ acts faithfully on a smooth curve
$C$ of genus 4, and this contradicts Lemma 3.7.

\th Proposition
\enonce Let $S$ be a del Pezzo surface admitting a
faithful action of $(\Z/3)^r$ with $r\ge 3$. Then $r=3$, $S$ is
isomorphic to the Fermat cubic $X_0^3+\ldots+X_3^3=0$, and
$(\Z/3)^3$ acts as the $3$\tx torsion
subgroup of the diagonal torus in $PGL_4(k)$.
\endth
{\it Proof} : The list 3.5 gives $\ell \ge 6$. On the other
hand (3.6) and Lemma 3.7 give $\ell \le 6$. Thus $S$
is a cubic surface in $\P^3$, and $G$ is a subgroup of $PGL_4(k)$. By
Lemma 3.1 we have $r= 3$, and
 there is a coordinate system $(X_0,\ldots ,X_3)$ on $\P^3$ such that
$G$ is the diagonal subgroup
$(\mug_3)^4/\mug_3$ of $PGL_4(k)$. 
\ind The surface  $S$ is defined by an element $F$of
$H^0(\P^3,{\cal O}_{\P^3}(3))$ which is semi-invariant with respect
to the action of $(\mug_3)^4$. Under this action the space
$H^0(\P^3,{\cal O}_{\P^3}(3))$ is the direct sum of the invariant
subspace spanned by
$X_0^3,\ldots ,X_3^3$ and of 16 1-dimensional subspaces spanned by
monomials, corresponding to 16 different characters of $(\mug_3)^4$.
Since $S$ is smooth, $F$ must be of the form $a_0X_0^3+\ldots
+a_3X_3^3$ with  $a_i\not= 0$ for each $i$, hence the result.\cqfd 

\th Proposition
\enonce Let $S$ be a del Pezzo surface and
$G$ a subgroup of $\aut(S)$ isomorphic to $(\Z/2)^r$ with $r\ge 4$ and
$\rk\Pic(S)^G=1$. Then $r=4$, $S$
is a quartic del Pezzo surface in
$\P^4$ with equations
$\sum_{i=0}^4X_i^2=\sum_{i=0}^4\lambda_iX_i^2=0$ for some distinct
elements $\lambda_0,\ldots \lambda_4$ of $k$, and $G$ is the  $2$\tx
torsion subgroup of the diagonal torus in $PGL_5(k)$.
\endth

{\it Proof} :  Once again the list 3.5 gives $\ell \ge 5$, and
(3.6) and Lemma 3.7 give $\ell \le 6$. 
\ind Suppose $S$ is a cubic surface. Then $G$ acts on the set of
lines in $S$; since 27 is odd, there must be one orbit with one
element, that is, one line stable under $G$. This contradicts the
assumption on $\Pic(S)^G$.
\ind Suppose $S$ is an intersection of two quadrics in
$\P^4$. Then
$G$ is a subgroup of $PGL_5(k)$; by Lemma 3.1 we have
 $r=4$, and  there is a coordinate
system $(X_0,\ldots ,X_4)$  on $\P^4$ such that
$G$ is the diagonal subgroup
$(\mug_2)^5/\mug_2$ of $PGL_5(k)$. 
\ind The representation of $(\mug_2)^5$ on $H^0(\P^4,{\cal
O}_{\P^4}(2))$ splits as
$$H^0(\P^4,{\cal O}_{\P^4}(2))=I\ \oplus\sum_{i<j} k\cdot  (X_iX_j)\ ,
$$where $I$ is the invariant subspace spanned by $X_0^2,\ldots
,X_4^2$. The 2-dimensional subspace of quadratic forms vanishing on
$S$ must be contained in $I$, since otherwise it would contain some
$X_iX_j$. After a change of coordinates we find the form given in the
Proposition.\cqfd 

\rem{Remark} The same method gives  the subgroups of type 
$(\Z/3)^2$ of the Cremona group: besides the two conjugacy classes of
subgroups of $PGL_3(k)$ described in (3.2), one gets the 
automorphisms groups of certain cubic surfaces and del Pezzo surfaces 
of degree 1. The subgroups of type
$(\Z/2)^3$ and
$(\Z/2)^2$ can also be described with analogous, but more tedious, 
methods.
\section{Conjugacy classes of 2-elementary subgroups}
\subsection It follows from the theorem that the $p$\tx elementary subgroups
of $\ca$ of maximal order form only one conjugacy class, except in the
case $p=2$. We are going to analyze the latter case. We assume ${\rm
char}(k)\not=2$ throughout this section.
\ind For $\sigma \in\ca$, we denote by
$NF(\sigma )$ the {\it normalized fixed locus} of $\sigma $, that is, the
 normalization of  the union of the non-rational curves in $\P^2$ fixed
by $\sigma $. The isomorphism class of $NF(\sigma )$ is an invariant of
the conjugacy class of $\sigma $, and this is the basic invariant we will use to
distinguish conjugacy classes.

\subsection Let $I$ be a finite subset of $k\smallsetminus \{2,-2\}$;
we associate to $I$ the subgroup $G_I$ of $\ca$ generated by the involutions
$$(x,y)\mapsto (-x,y)\quad,\quad (x,y)\mapsto ({1\over
x},y)\quad,\quad(x,y)\mapsto (x,-y)\quad,\quad(x,y)\mapsto (x,{P(x)\over
y})$$with $P(x)=\pprod_{a\in I}^{}(x^2+x^{-2}-a)$. By definition this gives the
 conjugacy
classes of subgroups  of type ${\rm c}_1)$.
 
\ind Consider the Galois $V_1$\tx covering\note{To distinguish the two copies
of ${\scriptstyle \P^1}$ which appear we indicate the coordinate by a
subscript.}
$\pi :\P^1_x\rightarrow
\P^1_u$ given by $u=x^2+x^{-2}$. The normalizer $N(V_1)$ of $V_1$ in
$PGL_2(k)$ is isomorphic to ${\goth S}_4$, and the quotient
$N(V_1)/V_1\cong {\goth S}_3$ acts on $\P^1_u$;  this   action 
preserves the branch locus 
$\{2,-2,\infty\}$ of $\pi $, and therefore identifies 
 ${\goth S}_3$ to the group of homographies 
of $\P^1_u$ preserving this subset. If $\gamma\in N(V_1)$, the birational
transformation $(x,y)\mapsto (\gamma(x), y)$ conjugates $G_I$ to
$G_{\gamma (I)}$. 
\th Proposition
\enonce The map $I\mapsto G_I$ induces a bijection between finite subsets of
$k\smallsetminus \{2,-2\}$ modulo the action of ${\goth S}_3$ and conjugacy classes
of subgroups of $\ca$ of type ${\rm c}_1)$.
\endth
{\it Proof} : We have already seen that the map is well-defined and surjective; it
remains to prove that we can recover from $G_I$ the finite subset $I\i k\smallsetminus
\{2,-2\}$ modulo ${\goth S}_3$.
\ind First observe that the case $I=\varnothing$ (which gives 
$G_I=V_1\times V_1\i\aut(\P^1\times \P^1)$) is the only one for which
$NF(\sigma )=\varnothing$ for all $\sigma \in G_I$. Fom now on we assume
$I\not=\varnothing$.
\ind For $\sigma \in G_I$,
we have $NF(\sigma )=\varnothing$ except for 
$\displaystyle \sigma _{\pm}:(x,y)\mapsto (x,\pm{ P(x)\over y})$;
$NF(\sigma _+)$ and $NF(\sigma _-)$ are isomorphic to the 
hyperelliptic curve $C_I$ with affine equation $y^2=P(x)$.  
The double covering $x:C_I\rightarrow \P^1_x$ is branched along $W=\pi ^{-1}
(I)$; the group $V_1$ acts freely on $W$.
\ind Conversely, starting from $G_I\i \ca$, we get the curve $C_I$, hence the
branch locus $W\i \P^1$ {\it up to the action of} $PGL_2(k)$. Moreover we
have an action of a  group  $G'\cong (\Z/2)^2$ on $\P^1$ preserving $W$:
indeed $G_I$ acts on 
$C=NF(\sigma _+)$, with $\sigma _+$ acting trivially and $\sigma _-$ acting as
the hyperelliptic involution. We take
  $G' =G_I/\langle \sigma _+,\sigma _-\rangle$.
\ind Choose a  basis $(e_1,e_2)$ of $G'$ over $\Z/2$. By Lemma 2.1
there is a coordinate $x$ on $\P^1$ such that $e_1$ acts by $x\mapsto -x$ and
$e_2$  by $x\mapsto x^{-1} $; it is unique up to the action of $V_1$. In
particular the map $\pi :x\mapsto x^2+x^{-2} $ and the subset
$I=\pi (W)$ of $k\smallsetminus \{2,-2\}$ depend only on the choice of the basis
$(e_1,e_2)$. A change of basis is realized by an element of $N(V_1)$, so the
class of $\pi (W)$ mod. ${\goth S}_3$ is well-defined. Thus we recover $I$ 
 mod. ${\goth S}_3$ from the subgroup $G_I$ of $\ca$.\cqfd
\smallskip 
\subsection We now consider  subgroups of type ${\rm c}_2)$.
We first recall some very classical facts about quartic del Pezzo
surfaces -- one possible reference is [H], lecture 22.
\ind  A quartic del Pezzo surface  $S\i\P^4$  is contained in a
pencil 
$(Q_\lambda )_{\lambda \in\P^1}$ of quadrics. There are exactly 5
singular quadrics $Q_{\lambda _0},\ldots Q_{\lambda _4}$ in this 
pencil; the map $S\mapsto \{\lambda _0,\ldots ,\lambda
_4\}$ is an isomorphism from the moduli space of quartic del Pezzo
surfaces onto the moduli space of 5-points subsets of $\P^1$ (modulo
the action of $PGL_2$). The  quadrics
$Q_{\lambda _0},\ldots Q_{\lambda _4}$ have rank
$4$; their singular points $p_0,\ldots p_4$ span $\P^4$. 
\ind  The group
$\aut(S)$ contains a normal, canonical subgroup $G_S$ isomorphic to 
$(\Z/2)^4$. Indeed for $0\le \ell \le 4$, there is a unique involution
$\sigma _\ell $ of $\P^4$ whose fixed locus consists of $p_\ell $ and the
hyperplane $H_\ell $ spanned by the points
$p_i$ for $i\not=\ell $; these involutions span the group $G_S$.
 In more concrete terms, choose the coordinates on $\P^1$ so that
$\lambda _0,\ldots ,\lambda _4\in k$. There exists a system of
coordinates on $\P^4$ such that the equations of $Q_{\infty}$ and
$Q_0$ are respectively
$$\sum_{i=0}^4 X_i^2=0\qquad \sum_{i=0}^4 \lambda _i\,X_i^2=0\
.$$Then  $G_S$ is  the $2$\tx torsion
subgroup of the diagonal torus in $PGL_5(k)$; the invo\-lution $\sigma
_\ell $ maps $(X_0,\ldots ,X_\ell ,\ldots ,X_4)$ to  $(X_0,\ldots ,-X_\ell
,\ldots ,X_4)$. As before we view
$G_S$ as a subgroup of $\ca$, well-defined up to conjugacy.
\th Proposition
\enonce {\rm a)} The subgroups of type ${\rm c}_1)$ are not conjugate to
those of type ${\rm c}_2)$.
\ind {\rm b)} Let $S$ and $S'$ be two quartic del Pezzo surfaces, with
$S$ general. If $G_{S'}$ is conjugate to $G_S$, $S'$ is isomorphic to $S$.
\endth
\ind In other words, the map $S\mapsto G_S$ from the moduli space of
quartic del Pezzo surfaces to the set of conjugacy classes of subgroups of
$\ca$  is generically injective. I~expect actually
this map to be injective, hence a bijection onto the set of conjugacy
classes of groups of type ${\rm c}_2)$; one way to prove it
would be by using classical invariant theory, but the computations look
rather frightening.\medskip 
{\it Proof of the Proposition} : Let $S$ be a quartic del Pezzo surface. The 5
elements $\sigma _\ell $ of $G_S$ fix the elliptic curve $H_\ell \cap S$, while
the other elements of $G_S$ have no fixed curve. On the other hand we
have seen that for a subgroup  $G\i\ca$ of type ${\rm c}_1)$, at most 2
elements of $G$ have a non-trivial normalized fixed locus. This is enough to
prove assertion a).
\subsection To prove b), we want  to recover the surface $S$ -- or, equivalently,
the 5-points subset $ \{\lambda _0,\ldots ,\lambda
_4\}$ -- from the data of the elliptic curves $H_\ell \cap S\ (0\le
\ell \le 4)$. The curve $H_\ell \cap S$ is given in $H_\ell \cong\P^3$
by the equations
$$\sum_{i\not=\ell }X_i^2= \sum_{i\not=\ell }\lambda
_i\,X_i^2=0\ ;$$it is isomorphic  to the double
covering of $\P^1$ branched along the points $\lambda _i$ for
$i\not=\ell $ ([H], Prop. 22.38). The  $j$\tx
invariant of this curve is
$\jmath([\lambda_0,\ldots ,\widehat{\lambda_\ell
},\ldots\lambda_4])$, where $[\ldots ]$ denotes the
cross-ratio and $\jmath$ is the
rational function
$\displaystyle \jmath(x )=2^8{(x ^2-x +1)^3\over
x ^2(x -1)^2}$. This function is invariant under the ``cross-ratio group",
that is, the subgroup of $PGL_2(k)$, isomorphic to ${\goth S}_3$,
generated by the involutions $z\mapsto z^{-1} $ and $z\mapsto 1-z$. In
fact the map $\jmath:\P^1\rightarrow \P^1$ is a Galois covering with
group ${\goth S}_3$. 

\subsection Let $V$ be the complement of the diagonals in $(\P^1)^5$. It
has two commuting actions of $PGL_2(k)$ and ${\goth S}_5$. We are
interested in the map $J:V\rightarrow (\P^1)^5$ which associates to a
5-tuple $(\lambda_0,\ldots ,\lambda_4)$ the 5-tuple $(j_0,\ldots
,j_4)$ given by $j_\ell =\jmath([\lambda_0,\ldots
,\widehat{\lambda_\ell },\ldots
\lambda_4])$. This map is ${\goth S}_5$\tx equivariant, and factors
through the quotient  $P:=V/PGL_2(k)$. Thus we have a commutative
diagram 
$$\diagram{P & \hfl{J}{} & (\P^1)^5 &\cr 
\vfl{}{} & & \vfl{}{}&\cr
P/{\goth S}_5 & \hfl{\bar J}{} & \Sy^5(\P^1)&\kern-8pt.}$$
The space $P/{\goth S}_5$ is the moduli space of 5-points sets in
$\P^1$, or equivalently of quartic del Pezzo surfaces;  the map
$\bar J$ associates to a quartic surface $S$ the 5 $j$\tx invariants of
the elliptic curves $H_\ell \cap S$. Our aim is to prove that  $\bar J$ is
generically injective; it suffices to prove that  $ J$ {\it is
generically injective}.

\subsection The quotient $P=V/PGL_2$
 is easy to describe: any 5-tuple is equivalent to one (and only one) of
the form $(\lambda,\mu,1,0,\infty)$. Thus we can identify $P$ with the
open subset of couples $(\lambda,\mu)$ in $\P^1\times \P^1$ with
$\lambda,\mu\notin\{0,1,\infty\}$ and $\lambda\not= \mu$. Then
the map $J$ is given by
$$J(\lambda,\mu)=\Bigl(\jmath(\mu)\,,\,\jmath(\lambda)\,,\,\jmath({\lambda
\over \mu })\,,\,\jmath({\lambda -1\over \mu -1})\,,\,\jmath({\lambda
(\mu -1)\over \mu (\lambda -1)})\Bigr)\ .$$ 
Since $\jmath^{-1} (\infty)=\{0,1,\infty\}$ we see that $J$ takes its
value in $(\A^1)^5=\A^5$.
 The strategy will be to first prove that $J:P\rightarrow \A^5$ is proper,
then that for some $p\in P$ the tangent map
$T_p(J)$ is injective and  $J^{-1}( J(p))=\{p\}$.
\subsection Let us view $J$ as a rational map $\P^1\times
\P^1\dasharrow(\P^1)^5$, obtained by composition of the rational map
$$f:\P^1\times
\P^1\dasharrow
(\P^1)^5\quad\hbox{given by}\quad 
f(\lambda,\mu)=\bigl(\mu\,,\,\lambda\,,\,{\lambda
\over \mu }\,,\,{\lambda -1\over \mu -1}\,,\,{\lambda
(\mu -1)\over \mu (\lambda -1)}\,\bigr)$$with the
morphism$\jmath^5:(\P^1)^5\rightarrow (\P^1)^5$. 
\ind The map $f$ is well-defined
except at the points
$(0,0),(1,1)\hbox{ and }(\infty,\infty)$, and extends to a morphism
$X\rightarrow (\P^1)^5$, where 
$X$ is the surface obtained by blowing up $\P^1\times
\P^1$ along these points. Thus $J$ extends  to $\bar J:X \rightarrow
(\P^1)^5$, which maps the three exceptional divisors in
$X$  into $\{\infty\}\times \{\infty\}\times
(\P^1)^3$. Therefore $\bar J$ maps $X\smallsetminus P$ to $(\P^1)^5\smallsetminus
(\A^1)^5$, so its restriction $J:P\rightarrow \A^5$ is proper.

\subsection We will now exhibit a point $p\in P$ such that $J^{-1}
(J(p))=\{p\}$. We choose an element $\zeta$ of $k$ with
$\zeta ^3=-1$; if
${\rm char}(k)\not=3$ we assume moreover $\zeta \not=-1$. 
We  take $p=(\alpha ,\zeta )$, with $\alpha $ general enough in $k$.
So  let
$(\lambda ,\mu )\not=(\alpha ,\zeta )$ in $P$ with $J(\lambda ,\mu
)=J(\alpha,\zeta  )$. Then $\mu $ belongs to the orbit $\{\zeta  ,\zeta ^{-1}
\}$  of 
$\zeta
$ under the cross-ratio group ${\goth S}_3$.
\ind Suppose $\mu  =\zeta $, and ${\rm char}(k)\not=3$. Then we have
$$\lambda  =\sigma (\alpha )\quad ,\quad  {\lambda  \over \zeta }=\tau
({\alpha
\over \zeta})\quad  \hbox{ for
some }\ \sigma ,\tau \not=1\ \hbox{ in }{\goth S}_3\ .$$
 This gives $\zeta ^{-1} \sigma (\alpha
)=\tau (\zeta ^{-1} \alpha )$; if $\alpha $ is general enough this implies
$h_\zeta^{-1}  \,\sigma \, h_\zeta  =\tau $ $ \in{\goth S}_3$,
where
$h_\zeta $ is the homothety $z\mapsto \zeta z$. Since ${\goth S}_3$
preserves $\{0,1,\infty\}$ we see that $\sigma $ must preserve
$\{0,\infty\}$, hence fix 1; the only  possibility is $\sigma 
(z)=z^{-1} $, which is excluded because $\zeta ^2\not=1$.
\ind Suppose $\mu =\zeta ^{-1} $ (with no assumption on ${\rm
char}(k)$). We have
$$\lambda  =\sigma (\alpha )\ ,\ {\zeta \lambda   }=\tau ({\alpha
\over \zeta})\ ,\ \zeta^{-1}  (1-\lambda  )=\xi (\zeta (1-\alpha ))\ \hbox{
for some }\ \sigma ,\tau ,\xi \ \hbox{ in }{\goth S}_3\ .$$
 The first two relations give $\zeta  \sigma (\alpha
)=\tau (\zeta ^{-1} \alpha )$, so that  $h_\zeta \,\sigma \, h_\zeta 
=\tau  \in{\goth S}_3$. As before the only possibility for $\sigma $ is
$\sigma (z)=z^{-1} $, so that $\lambda =\alpha ^{-1} $. Now the third
relation becomes $\zeta ^{-1} (1-\alpha ^{-1} )=\xi (\zeta (1-\alpha ))$;
this cannot hold for $\alpha $ generic because it does not hold for
$\alpha =\infty$.
\subsection It remains to check that the tangent map to $J$ at
$p=(\alpha ,\zeta )$ is injective for $\alpha $ general enough.
 The Jacobian determinant$$  \det\pmatrix{{\partial \over \partial
\lambda }\jmath({\lambda \over\mu }) &{\partial \over \partial
\mu  } \jmath({\lambda \over \mu})\cr
{\partial \over \partial
\lambda }\jmath({\lambda -1\over\mu -1 }) &{\partial \over \partial
\mu  } \jmath({\lambda  -1\over \mu -1}) }=\jmath'({\lambda \over\mu
})\,\jmath'({\lambda -1\over\mu -1}){\mu -\lambda \over \mu ^2(\mu
-1)^2}$$is nonzero if ${\lambda\over\mu}$ and ${\lambda
-1\over\mu -1}$ do not belong to $\{\zeta  ,\zeta^2,-1,2,{1\over
2}\}$, which is certainly the case when $\alpha $ is general enough. This
finishes the proof of the Proposition.\cqfd

\vskip2truecm
\centerline{ REFERENCES} \vglue15pt\baselineskip12.8pt
\def\num#1{\smallskip\item{\hbox to\parindent{\enskip [#1]\hfill}}}
\parindent=1,15cm

\num{B-B} L. {\pc BAYLE}, A. {\pc BEAUVILLE}: {\sl Birational involutions
of} $\Bbb{P}^2$. Kodaira's issue, Asian J. Math. {\bf 4} (2000) 11--17. 

\num{B-Bl} A. {\pc BEAUVILLE}, J. {\pc BLANC}: {\sl On Cremona
transformations of  prime order}.  C. R. Math. Acad. Sci. Paris {\bf 339}
(2004), no. 4, 257--259.
\num{B}  N. {\pc BOURBAKI}: {\sl Groupes et alg\`ebres de Lie}, Chap.\
VI. Hermann, Paris (1968).

\num{Bo} A. {\pc BOREL}: {\sl Sous-groupes commutatifs et torsion des
groupes de Lie compacts connexes}. T\^ohoku Math. J. (2) {\bf 13}
(1961), 216--240.
\num{D1} M. {\pc DEMAZURE}: {\sl Sous-groupes alg\'ebriques de rang
maximum du groupe de Cremona}.  Ann. Sci. \'Ecole Norm. Sup.
(4) {\bf 3} (1970), 507--588.
\num{D2} M. {\pc DEMAZURE}: {\sl Surfaces de del Pezzo}. S\'eminaire
sur les Singularit\'es des Surfaces. LNM 
{\bf 777}, 21--69. Springer, Berlin (1980).
\num{Do} I. {\pc DOLGACHEV}: {\sl Weyl groups and Cremona
transformations}.  Singularities I, 283--294, 
Proc. Sympos. Pure Math. {\bf 40}, 
AMS, Providence (1983).
\num{dF} T. {\pc {}DE} {\pc FERNEX}: {\sl On planar Cremona maps of
prime order}. Nagoya Math. J. {\bf 174} (2004), 1--28. 
\num{dF-E} T. {\pc {}DE} {\pc FERNEX},  L. {\pc EIN}: {\sl Resolution of
indeterminacy of pairs}. Algebraic geometry, 165--177, de Gruyter, 
Berlin (2002). 
\num{H} J. {\pc HARRIS}: {\sl Algebraic geometry,
a first course}.  Graduate Texts in
Mathematics, {\bf 133}.  Springer-Verlag, New York (1995).
\num{K} S. {\pc KANTOR}: {\sl  Theorie der endlichen Gruppen von
eindeutigen Transformationen in der Ebene}. Mayer \& M\"uller, Berlin
(1895).
\num{M} Y. {\pc MANIN}: {\sl Rational surfaces over perfect fields}, II. Math.
USSR -- Sbornik {\bf 1} (1967), 141--168.
\num{W} A. {\pc WIMAN}: {\sl Zur Theorie der endlichen Gruppen von
birationalen Transformationen in der Ebene}. Math. Ann. {\bf 48} (1896),
195--240.
\num{Z} D.-Q. {\pc ZHANG}: {\sl Automorphisms of finite order on
rational surfaces}.  J. Algebra {\bf 238}
(2001), no. 2, 560--589. 
\vskip1cm
\def\pc#1{\eightrm#1\sixrm}
\hfill\vtop{\eightrm\hbox to 5cm{\hfill Arnaud {\pc BEAUVILLE}\hfill}
 \hbox to 5cm{\hfill Institut Universitaire de France\hfill}\vskip-2pt
\hbox to 5cm{\hfill \&\hfill}\vskip-2pt
 \hbox to 5cm{\hfill Laboratoire J.-A. Dieudonn\'e\hfill}
 \hbox to 5cm{\sixrm\hfill UMR 6621 du CNRS\hfill}
\hbox to 5cm{\hfill {\pc UNIVERSIT\'E DE}  {\pc NICE}\hfill}
\hbox to 5cm{\hfill  Parc Valrose\hfill}
\hbox to 5cm{\hfill F-06108 {\pc NICE} Cedex 02\hfill}}
\end